\documentclass{amsart}
\usepackage{amssymb,amscd,amsmath,hyperref}
\usepackage{graphicx}
\theoremstyle{plain}
\newtheorem{theorem}[subsection]{Theorem}
\newtheorem{cor}[subsection]{Corollary}
\newtheorem{lem}[subsection]{Lemma}
\newtheorem{prop}[subsection]{Proposition}
\theoremstyle{definition}
\newtheorem{defn}[subsection]{Definition}
\theoremstyle{definition}
\newtheorem{rem}[subsection]{Remark}
\theoremstyle{definition}

\theoremstyle{definition}

\DeclareMathOperator{\N}{N}
\DeclareMathOperator{\Ne}{N}

\DeclareMathOperator{\interior}{int}
\DeclareMathOperator{\ver}{ver}

\DeclareMathOperator{\supp}{supp}
\DeclareMathOperator{\Hom}{Hom}
\begin{document}
\title[Alexander and Thurston norms of graph links]%
{Alexander and Thurston norms of graph links}

\author[David G.~ Long]{David G. Long}
\address{Department of Mathematics, Northeastern University,
Boston, MA, 02115 }
\email{\href{mailto:dlong2002us@yahoo.com}{dlong2002us@yahoo.com}}


\begin{abstract}
We show that the Alexander and Thurston norms are the same for all irreducible Eisenbud-Neumann graph links in homology $3$-spheres.  These are the links obtained by splicing Seifert links in homology $3$-spheres together along tori. By combining this result with previous results, we prove that the two norms coincide for all links in $S^3$ if either of the following two conditions are met; the link is a graph link, so that the JSJ decomposition of its complement in $S^3$ is made up of pieces which are all Seifert-fibered, or the link is alternating and not a $(2,n)$-torus link, so that the JSJ decomposition of its complement in $S^3$ is made up of pieces which are all hyperbolic. We use the E-N obstructions to fibrations for graph links together with the Thurston cone theorem on link fibrations to deduce that every facet of the reduced Thurston norm unit ball of a graph link is a fibered facet. \end{abstract}
\maketitle
\section{Introduction}
\subsection{The Alexander norm (A-norm)}  A semi-norm defined on the first cohomology group $H^1(M;\mathbb{Z})$ of a connected, compact, orientable 3-manifold $M$, whose boundary (if any) is a union of tori, was introduced by McMullen \cite{McM} in the late 1990s. It is directly determined by the multivariable Alexander polynomial $\triangle$ of the $3$-manifold $M$ and is called the \textit{Alexander norm} of  $M$. We shall adopt the notation A-norm for the Alexander norm. 

If we write a multivariable Alexander polynomial $\triangle$ in multi-index notation, then $\triangle = \sum_{\alpha \in \supp(\triangle)} c_{\alpha} \textbf{t}^{\alpha} $ with $\textbf{t}^{\alpha} = t_1^{\alpha_1} \dots t_n^{\alpha_n}$ and $\supp(\triangle)$ denoting the support of the polynomial $\triangle$; the set of all $\alpha$ labeling non-zero constants $c_{\alpha}$.  Let $\alpha$ and $\beta$ be the exponents of any two arbitrarily chosen terms of $\triangle$.  They are elements of $H_1(M; \mathbb{Z})$.   Let $\phi$ be a cohomology class in the space dual to the space containing the exponents so that $\phi \in H^1(M;\mathbb{Z}) =\mathbb{Z}^{b_1}$ where $b_1$ is the first Betti number of $M$.  The A-norm $\| \phi \|_A$ of $\phi$ is the supremum of $ \phi(\alpha - \beta)$ taken over all the exponents $\alpha, \beta$ in the support of $\triangle$. 

\subsection{The Thurston norm (T-norm)}The A-norm is closely related to a semi-norm for compact, oriented $3$-manifolds introduced in 1986 in \cite{Thu}, called the \textit{Thurston norm}.  We shall adopt the notation T-norm for the Thurston norm. In the compact, oriented 3-manifold $M$ with a boundary $\partial M$ which may be empty, every homology class $ z \in H_2(M,\partial M;\mathbb{Z})$ can be represented as $[S]$ where $S$ is a properly embedded, oriented surface. By Poincar\'e duality the relative second homology group is isomorphic to the first cohomology group, $H_2(M,\partial M;\mathbb{Z}) \cong H^1(M,\mathbb{Z})$. The $T$-norm  $\| \phi \|_T$ is determined by the Euler characteristic $\chi$ of the surface $S$ representing the cohomology class $\phi$ through this duality isomorphism.

\subsection{Eisenbud-Neumann (E-N) graph links} An E-N graph link (see \cite{Eis, Neu}) is a link $L=L_1 \cup \dots \cup L_r$ in a  homology $3$-sphere $\Sigma$, which can be built up by splicing together Seifert links in homology $3$-spheres. We use the notation $\mathbf{L} = (\Sigma, L)$ to indicate not only the link $L$ but also the homology 3-sphere $\Sigma$ it lies in. To each graph link $\mathbf{L}$ there is a unique minimal graph. It is shown in \cite{Eis} that the graph directly determines the Alexander polynomial and Thurston norm of  $\mathbf{L}$.

\subsection{Main theorem: The coincidence of Alexander and Thurston norms for graph links}
\begin{theorem} The Alexander and Thurston norms coincide for all irreducible graph links with two or more components.  \end{theorem}
We prove this theorem by directly calculating the A-norm of graph link $\mathbf{L}$ from the expression for its Alexander polynomial, and comparing it to the expression for its T-norm.   We derive three new results which are used in the proof. Each of these results involves a sum over the $p$ splice components of the graph of $\textbf{L}$;\begin{enumerate}\item a \textit{T-norm decomposition formula}, 
\item  the Newton polyhedron of the Alexander polynomial of $\mathbf{L}$ as a Minkowski sum, and 
\item an \textit{A-norm decomposition formula}. \end{enumerate}
In the case of a graph knot, which is a graph link with one component, we show as a corollary that the two norms satisfy $\| \phi \|_A - \| \phi \|_T = | \phi |$. 

Without actually calculating the two norms themselves, McMullen proved in \cite{McM} that for a compact, connected, oriented $3$-manifold $M$ whose boundary if any is a union of tori, the two norms must coincide for all fibered cohomology classes which are primitive. A primitive cohomology class $\phi$ is one whose image is $\mathbb{Z}$ so that $\phi(H_1(M;\mathbb{Z})) = \mathbb{Z}$. Hence our result that the two norms coincide for graph links is new with regard to non-primitive cohomology classes and the non-fibered cohomology classes. 

In \cite{McM}, McMullen states that the two norms, Alexander and Thurston, coincide for many 3-manifolds but not all. Examples of links for which the two norms do not coincide are given by a two-component link found by Dunfield \cite{Dun} and the link $9_{21}^3$ (See \cite{McM}).
\subsection{Generalization of the main theorem}As an application of our main result, we combine it with previous results to obtain a more general theorem in the context of the JSJ (Jaco, Shalen and Johannson) decomposition of a link in $S^3$. According to the JSJ decomposition, the link complement $X$ can be decomposed into irreducible pieces of two types, Seifert-fibered and hyperbolic, by de-splicing the pieces together along tori. 
\begin{theorem} Let $L$ be a link in $S^3$.  Let the JSJ decomposition of the link complement $S^3 \setminus L$ consist of $k$ irreducible pieces $\mathcal{N}_i$ so that \begin{equation} \nonumber S^3 \setminus L = \mathcal{N}_1 \cup_T \mathcal{N}_2 \cup_T \dots \cup_T \mathcal{N}_k. \end{equation}
Then the Thurston and Alexander norms of $L$ coincide if either of the following two conditions are met:
\begin{enumerate}
\item $L$ is a graph link so that all the $\mathcal{N}_i$ are Seifert-fibered.
\item $L$ is an alternating link which is not a $(2,n)$-torus link, so that all the $\mathcal{N}_i$ are hyperbolic. \end{enumerate}\end{theorem}
\subsection{Characteristic hyperplanes of the Thurston norm ball}
\subsubsection{T-norm unit ball and fibrations}The T-norm can be extended by convexity from integer-valued classes to real-valued classes to determine a T-norm unit ball. The cohomology classes of the T-norm unit ball can be separated into two types; fibered and non-fibered. An integer-valued cohomology class $\phi$ is fibered if the surface $S$ which represents it via Poincar\'e duality is also the fiber $F$ of a fibration of 3-manifold $M$ over the circle $S^1$. A real-valued cohomology class $\phi$ is fibered if it lies on a ray through the origin whose lattice points are integer-valued fibered cohomology classes. Otherwise $\phi$ is a non-fibered cohomology class. We do not require the fiber to be connected. 
\subsubsection{Thurston cone theorem} The Thurston norm unit ball may not be a bounded set. However, we show in Long \cite{Lon} that by introducing \textit{essential coordinates} we can define a \textit{reduced T-norm unit ball} which is a polyhedron of the same dimension as the Newton polyhedron of the Alexander polynomial of 3-manifold $M$. The \textit{Thurston cone theorem} (see \cite{Thu}) states that the set of fibered cohomology classes is some union of the cones pointed at the origin (minus the origin) through the interiors of the facets (top-dimensional faces) of the reduced T-norm unit ball.

\subsubsection{Fibration obstruction criterion}  Using necessary and sufficient conditions for a cohomology class of a graph link $\mathbf{L}$ to be a fibered class from \cite{Eis}, we show that the set of all non-fibered classes of $\mathbf{L}$ is made up of a set of hyperplanes we call the \textit{characteristic hyperplanes}. There is one characteristic hyperplane for each splice component of the graph of $\mathbf{L}$. The $i$th characteristic hyperplane is made up of those classes whose T-norm for the $i$th splice component in the Thurston norm decomposition formula is zero. 
\subsubsection{The fibered facets of the reduced Thurston norm unit ball} We use the fibration obstruction criterion and the Thurston cone theorem to prove the following new theorem on graph link fibrations.\begin{theorem}  Every facet $\mathcal{F}$ of the reduced T-norm unit ball $\mathcal{\tilde{B}}_T$ for the irreducible graph link $\textbf{L}$ is a fibered facet. \end{theorem}
This theorem implies that the boundary of each facet, which is made up of the lower dimensional faces, must be contained in the non-fibered set of cohomology classes. 
\subsection{Sample calculation} We conclude by applying our results to a sample graph link $\mathbf{L}_{EN}$ from Eisenbud and Neumann \cite{Eis}.  We find the Alexander polynomial, Thurston norm, reduced T-norm unit ball and characteristic hyperplanes for this link and determine the intersection of the characteristic hyperplanes with the reduced T-norm unit ball.
\section{Alexander and Thurston norms} 
\subsection{Alexander norm} Let $\triangle$ be the Alexander polynomial of a compact, connected, orientable 3-manifold $M$ with first Betti number $b_1$, whose boundary (if any) is a union of tori. This polynomial can be expressed using multi-index notation as follows: 
\begin{eqnarray}\label{E:pol}  \triangle(t_1, \dots ,t_{b_1}) & =& \sum c_{\alpha_1 \alpha_2 \dots \alpha_{b_1}} t_1^{\alpha_1}t_2^{\alpha_2} \dots t_{b_1}^{\alpha_{b_1}}\\ \nonumber
& =& \sum_{\alpha \in \supp(\triangle)} c_{\alpha} \textbf{t}^{\alpha} 
  \end{eqnarray} 
$\mbox{where } \supp(\triangle) = \lbrace \alpha \colon c_{\alpha} \neq 0 \rbrace.$ The A-norm is a semi-norm defined on the first cohomology group $H^1(M;\mathbb{Z})$ of a 3-manifold $M$ directly determined by the  Alexander polynomial $\triangle$ of $M$. The Alexander polynomial $\triangle$ is a Laurent polynomial; $\triangle \in  \mathbb{Z}[t_1^{\pm 1}, \dots , t_{b_1}^{\pm 1}]$. Hence $\triangle$ can be expressed as a finite sum as in Equation \eqref{E:pol}. Then $\alpha \in \mathbb{Z}^{b_1}, \forall \alpha \in \supp(\triangle)$.
Further let $\phi \in (\mathbb{Z}^{b_1})^* = \Hom_{\mathbb{Z}}(\mathbb{Z}^{b_1}, \mathbb{Z})=H^1(M; \mathbb{Z}) \cong \mathbb{Z}^{b_1}$ be an element of the dual vector space $(\mathbb{Z}^{b_1})^*$ to the space $\mathbb{Z}^{b_1}$ in which each $\alpha$ lies in. Then there is a semi-norm, which we shall call the A-norm, defined as follows.
\begin{defn}\label{D:fn} The A-norm of $\phi \in H^1(M; \mathbb{Z})$ for the Alexander polynomial $\triangle$ of $M$ is  \begin{eqnarray}\label{E:fn}  && \| \hspace{2mm} \|_A \colon (\mathbb{Z}^{b_1})^* \to \mathbb{R}^+ \cup 0. \\ \nonumber
&&  \| \phi \|_A = \sup_{\alpha,\beta \in \supp(\triangle)} \phi(\alpha - \beta). \end{eqnarray}\end{defn}
\begin{rem} The definition of the A-norm above can be extended to real-valued cohomology classes $\phi \in H^1(M;\mathbb{R})$ =$(\mathbb{R}^{b_1})^*$. \end{rem}
Since the A-norm is completely determined by the Alexander polynomial $\triangle$ we also use the notation $\| \phi \|_A: = \| \phi \|_{\triangle}$.  In Long \cite{Lon}, we note that the Alexander norm is the special case of a norm determined by Laurent polynomials in general for which the polynomial is an Alexander polynomial. We call this generalized norm the Laurent norm since each Laurent polynomial $f$ with integer coefficients determines a norm $\| \hspace{2mm}\|_f$; the Laurent norm for $f$. 
 
\subsection{Thurston norm}The T-norm for compact, oriented 3-manifolds (with or without a boundary $\partial M$) was first defined on the second relative homology group $H_2(M,\partial M; \mathbb{Z})$ of $M$ in Thurston \cite{Thu}.  Any class of this relative second homology group can be represented by a compact, oriented two-dimensional surface $S$ in $M$. Each such surface has an integer-valued Euler characteristic $\chi(S)$ which can be used to define a norm on $H_2(M,\partial M; \mathbb{Z})$. By Poincar\'{e} duality each class of the group $H_2(M, \partial M; \mathbb{Z})$ determines a class of the first cohomology group $H^1(M;\mathbb{Z})$ of $M$;  $H_2(M,\partial M; \mathbb{Z}) \cong H^1(M;\mathbb{Z})$. Due to this equivalence, the T-norm can also be defined as a norm on $H^1(M;\mathbb{Z})$. We use the formulation of the T-norm as a norm on $H^1(M;\mathbb{Z})$.  In Dunfield \cite{Dun}, this version of the T-norm is defined as follows: 
\begin{defn}[Dunfield, \cite{Dun}]For a compact, connected surface $S$, let $\chi^-(S) = | \chi(S)|$ if $\chi(S) \leq 0$ and 0 otherwise. For a surface with multiple connected components $S_1, S_2, \dots ,S_n$ let $\chi^-(S)$ be the sum of the $\chi(S_i)$. Then the \textit{T-norm} of an integer-valued class $\phi \in H^1(M;\mathbb{Z}) \cong H_2(M,\partial M; \mathbb{Z})$ is  \begin{equation} \nonumber || \phi ||_T = \left\{\begin{array}{l}\inf  \chi^-(S)\mid S \mbox{ is a properly embedded oriented} \\
 \mbox{surface that is dual to } \phi  \end{array}\right\} .\end{equation}  \end{defn}
 
\begin{rem} The surface $S$ of this definition may or may not be connected. It is shown in \cite{Thu} that if cohomology class $\tilde{\phi}$ satisfies $\tilde{\phi} = d \cdot \phi$ for some integer $d \geq 1$, then $\| \tilde{\phi} \|_T = 
d \cdot \| \phi \|_T$ because $\tilde{\phi}$ represents $d$ disjoint surfaces each with T-norm $\| \phi \|_T$. The T-norm is additive for disjoint surfaces. \end{rem}
Also in \cite{Thu}, the T-norm is extended using convexity from integer-valued classes $\phi \in (\mathbb{Z}^{b_1})^*$ to real-valued classes $\phi \in (\mathbb{R}^{b_1})^*$ to determine a convex set in $(\mathbb{R}^{b_1})^*$ called the T-norm unit ball $\mathcal{B}_T$.  Thus we consider the T-norm as a semi-norm on $H^1(M;\mathbb{R})$, the first cohomology group of $M$ with real coefficients, rather than $H^1(M;\mathbb{Z})$.  

\section{E-N graph links}
\subsection{Link splicing} We use the notation that $\mathbf{L}=(\Sigma,L)$ denotes the link $L$ and its ambient space $\Sigma$ which is a homology $3$-sphere. Given two links, $L$ in $\Sigma$ and $L'$ in $\Sigma'$, with $r$ and $s$ components respectively we may form the link $L''$ in $\Sigma''$ with $r+s-1$ components by selecting link components, $S$ of $L$ and $S'$ of $L'$ and splicing the two links together along $S$ and $S'$.  To do this, first assume that $S$ and $S'$ have tubular neighborhoods with meridian and longitude $(m,l)$ and $(m',l')$ respectively. To splice $L$ to $L'$ along $S$ and $S'$ we attach $m$ to $l'$ and $l$ to $m'$. \begin{eqnarray}\nonumber  & & m \to l',\\ \nonumber & & l \to m'. \end{eqnarray} 
The splice is a homeomorphism of a tubular neighborhood of $S$ which is a solid torus, $D^2 \times S^1$, to a tubular neighborhood of $S'$, also a solid torus, $S^1 \times D^2$, along a boundary torus $S^1 \times S^1$.  The action can be represented as $D^2 \times S^1 \cup_{S^1 \times S^1} S^1 \times D^2$. This union of two solid tori across a torus induces the appropriate map in homology connecting each meridian with a longitude. 

\subsection{Seifert links} 
\subsubsection{Seifert-fibered spaces}  To construct a Seifert-fibered space, we start with a solid torus $T$ and remove $n$ parallel tori to obtain the space $T_n$. The fundamental group of this space $T_n$ is given by \begin{equation}\nonumber \pi_1 T_n = \langle h,y_0,y_1, \dots ,y_n \mid [h,y_i]=1, (0 \leq i \leq n), y_0 \cdot y_1 \cdots y_n=1 \rangle. \end{equation} This is an $S^1$-fibration of $T_n$ with base space a disc with $n$ points removed and circle fibers.  $h$ and $y_0$ represent the longitude and meridian of the ambient solid torus $T$.  Each $y_i, i \neq 0$, represents a meridian of $i$th torus determined by the $i$th hole of the disc. The commutation relations say that each generator $y_i, i=0, \dots,n $ commutes with the longitude $h$ of $T$ which means the fibration has trivial monodromy. Next we attach $n+1$ solid tori, $V_i$, to $T_n$ along the boundary tori in a special way using attaching maps that are homeomorphisms. \begin{equation}\nonumber M(e;(\alpha_1,\beta_1), \dots ,(\alpha_n,\beta_n)) = T_n \cup  \bigcup_{i=0}^{n} V_i .\end{equation} Each meridian generator, $m_i$, of the solid torus, $V_i$, is mapped  to $\alpha_i y_i + \beta_i h$ for $i \neq 0$.  For the solid torus $V_0$ the meridian generator, $m_0$  is mapped to $y_0 + e h $ where $e$ is an integer: \begin{eqnarray}\nonumber  && m_i = \alpha_i y_i + \beta_i h , i= 1, \dots , n .\\ \nonumber && m_0 = s_0 + e h .  \end{eqnarray}  
The space we obtain by this procedure is called the Seifert-fibered space \\ \textit{$M(e;(\alpha_1,\beta_1), \dots ,(\alpha_n,\beta_n))$}. It has $n$ exceptional fibers of type $(\alpha_i,\beta_i), i=1, \dots ,n)$ where $\alpha_i$ and $\beta_i$ denote the number of times the $i$th exceptional fiber wraps around the $i$th torus longitudinally and meridianally respectively. The exceptional fibers are the core circles of the solid tori after we have attached them; if $V_i = D^2 \times S^1$, then the core circle of $V_i$ which becomes the $ith$ exceptional fiber is $0 \times S^1$.
This space has fundamental group with presentation \begin{eqnarray} \nonumber \pi_1(M(e;(\alpha_1,\beta_1), \dots ,(\alpha_n,\beta_n))) &= & \langle  h,y_1, \dots , y_n \mid y_i^{\alpha_i} h^{\beta_i}=1, [h,y_i,]=1, \\ \nonumber &&
 (0 \leq i \leq n), y_1 \cdots y_n h^{-e} =1\rangle . \end{eqnarray} It is equivalent to the fundamental group of $T_n$ with the one additional relation added for each of the solid tori $V_i, i=0, \dots n$ attached. This construction of Seifert-fibered spaces is taken directly from Zieschang \cite{Ziec}. 
\subsubsection{Seifert link: }We obtain an $r$-component \textit{Siefert link} from the Seifert-fibered space $M(e;(\alpha_1,\beta_1), \dots , (\alpha_n, \beta_n))$ by removing a tubular neighborhood from each of the first $r$ exceptional fibers. Each component $S_i$, with $i=1, \dots ,r$, has a complement in the ambient Seifert-fibered space which is equivalent to a torus knot, labeled by $(\alpha_i, \beta_i)$. We leave the remaining $n-r$ components alone so that the complement of link has $n-r$ singular fibers. This link is denoted $L=S_1 \cup \dots \cup S_r$. The ambient Seifert-fibered space and $L$ together form the pair $\mathbf{L}=(M(e;(\alpha_1,\beta), \dots , (\alpha_n, \beta_n)),S_1 \cup \dots \cup S_r)$. Thus $L$ denotes the link and $\textbf{L}$ denotes both the ambient space $\Sigma=(M(e;(\alpha_1,\beta), \dots , (\alpha_n, \beta_n))$ the link lies in and the link $L$ itself; $\mathbf{L}= (\Sigma,L)$. The following is an equivalent but more concise definition of a Seifert link. 
\begin{defn}[Eisenbud and Neumann, \cite{Eis}, pg. 24] Let $L$ be a link in a $3$-manifold $M$ and let the interior of a closed tubular neighborhood $N(L)$ of $L$ be denoted $\interior N(L)$. Then $\textbf{L} = (M,L)$ is a \textit{Seifert link} if the link exterior, which is $M \setminus \interior N(L)$, of $L$ in $M$ is a Seifert-fibered space. \end{defn} 
  
\begin{rem} If the Seifert link $L$ is in the 3-sphere $S^3$, we drop the bold-faced notation so that $\textbf{L}= (S^3, L): = L$ in this case. \end{rem}
\subsection{Seifert-fibered homology 3-spheres}If we require that the Seifert-fibered space has the homology of the $3$-sphere we must have the additional relations that 
\begin{equation}\nonumber  \sum_{i=1}^n \beta_i \alpha_1 \dots \hat{\alpha_i} \dots \alpha_n = 1\mbox{ and }
  e=0. \end{equation} 
The notation $\hat{\alpha_i}$ means to remove the component $\alpha_i$ from the equation. In this case it is not necessary to include the coefficients $\beta_i, i=1, \dots ,n$ in uniquely labeling the Seifert-fibered space so that if the space is in a homology 3-sphere we use the notation $\Sigma(\alpha_1, \dots ,\alpha_n)$ and call it an \textit{unoriented} Seifert-fibered homology 3-sphere of type $(\alpha_1, \dots , \alpha_n)$. We use the notation $\Sigma(\epsilon;\alpha_1, \dots , \alpha_n)$ with $\epsilon=\pm 1$ to indicate the two possible orientations of the Seifert-fibered homology 3-sphere. Thus we use \textit{oriented} Seifert-fibered homology 3-spheres. The ambient oriented Seifert-fibered homology 3-sphere   together with link $L=S_1 \cup \dots \cup S_r$ form the pair denoted  $\mathbf{L}=(\Sigma(\epsilon;\alpha_1, \dots , \alpha_n, ),S_1 \cup \dots \cup S_r)$. This is the notation we shall use for all Seifert links which we assume are in oriented Seifert-fibered homology 3-spheres. The Seifert link $\mathbf{L}$ as defined above is an $r$-component link with first homology group given by $\mathbb{Z}^r$ as would be the case if it was instead in the $3$-sphere. To a given unoriented Seifert-fibered homology sphere, there exists a unique unordered $n$-tuple of coprime integers $(\alpha_1, \dots ,\alpha_n)$  with $\alpha_i \geq 2 , \forall i$.

\subsection{Graphs of Seifert links and splice diagrams} We can represent a Seifert link in a Seifert-fibered homology 3-sphere as a graph with the following components.
\begin{enumerate}
\item Boundary vertices: $\-----\bullet$ 

The boundary vertex corresponds to the solid torus which is a neighborhood of an exceptional fiber labeled by $\alpha_i$. 
\item Arrowhead vertices:  $\---- \longrightarrow$

The arrowhead vertex corresponds to one of the link components where a tubular neighborhood of an exceptional fiber has been removed. It also has the label $\alpha_i$.
\item Nodes: $\oplus$ and $\ominus$ 

A node corresponds to a Seifert manifold embedded in a link exterior;  the edges incident to a node correspond to the boundary components of the Seifert manifold, or for edges that lead to boundary vertices or arrowhead vertices to boundaries of tubular neighborhoods of fibers.  Within each node we insert $+$ and $-$ corresponding to the two possible orientations of the Seifert manifold. Each node must have at least three edges incident on it.
\end{enumerate}

We can now form a splice diagram (graph) by connecting two arrowhead vertices of the graphs of a pair of Seifert-fibered links. An example of a splice diagram is shown in Figure \ref{F:S2}. 
A splice diagram determines an E-N graph link which is defined as follows:
\begin{defn} An \textit{E-N graph link} consists of either a Seifert link in a homology 3-sphere or the splice of two or more Seifert links in homology 3-spheres. \end{defn}
In the rest of this article, we assume that the graph of graph link $\textbf{L}$ has $r$ arrowhead vertices labeled $v_1, \dots ,v_r$, $p$ nodes labeled $v_{r+1}, \dots ,v_{r+p}$ and $q$ boundary vertices labeled $v_{r+p+1} ,\dots ,v_{r+p+q}$. We also assume that the graph link is \textit{irreducible}; it can not be expressed as a disjoint sum. In addition, by the term graph link we mean an E-N graph link. 
\subsection{Alexander polynomial of a graph link}
The Alexander polynomial of the graph link $\textbf{L}$ is directly determined by its graph $\Gamma$.
\begin{theorem}[Eisenbud and Neumann \cite{Eis}, Theorem 12.1] Assume $\Gamma$ (the graph of $\textbf{L}$) is connected. Then the Alexander polynomial of the graph link $\textbf{L}$ is \begin{equation}\label{E:Al1} \triangle^\textbf{L}(t_1, \dots ,t_r)=\prod_{j=r+1}^{r+p+q}(t_1^{l_{1j}} t_2^{l_{2j}} \dots t_r^{l_{rj}} -1)^{\delta_j-2}.\end{equation}If $r =1$, so that the link is a knot, the formula is \begin{equation}\label{E:Al1k} \triangle^\textbf{L}(t_1) = (t_1 - 1)\prod_{j=2}^{p+q+1} (t_1^{l_{1j}} - 1)^{\delta_j - 2}. \end{equation}  Any terms of the form $(t_1^0 \dots t_r^0-1)^d$ ($d$ denotes an integer) which may occur on the right-hand side of these two equations should be formally canceled against each other before being set equal to zero. \end{theorem}
The Alexander polynomial has the same number of variables as arrowhead vertices and the product is over all vertices which are not arrowhead vertices; the nodes and boundary vertices. $\delta_j$ indicates the number of edges incident of each vertex.
The linking numbers $l_{ij}$ can be obtained directly from the graph.  For any two distinct vertices $v_i$ and $v_j$ of a graph $\Gamma$, let $\sigma_{ij}$ be the simple path in $\Gamma$ joining $v_i$ to $v_j$, including $v_i$ and $v_j$.  Then we have that 
\begin{equation} \label{E:ln}  l_{ij}   = \left\{ \begin{array}{l}  \mbox{(product of all signs of nodes on } \sigma_{ij} )   \cdot \mbox{(product of all edge weights}\\\mbox{adjacent to these nodes but not on } \sigma_{ij}). \end{array} \right\}\end{equation}
\subsection{T-norm of graph a link}
The T-norm of the graph link $\textbf{L}$ is also directly determined by the graph of $\textbf{L}$.
\begin{theorem}[Eisenbud-Neumann \cite{Eis}, Thm.\ 11.1] \label{T:TN} The T-norm of $\phi$ for the irreducible graph link $\textbf{L}$, which is not the unknot in $S^3$, is \begin{equation}\label{E:TN}  \|\phi \|_T = \sum_{i=r+1}^{r+p+q} (\delta_i - 2) \big{|}. \sum_{j=1}^r \phi_j  l_{ji}\big{|} 
 \end{equation} \end{theorem}
\begin{rem} For a graph link with $r$ components, let $d$ be the greatest common divisor of the $r$ components of the vector $\phi \in \mathbb{Z}^r$.   Up to homeomorphism, the surface $S$ that the cohomology class $\phi$ represents is the disjoint union of $d$ identical, connected, compact, oriented, two-dimensional surfaces each having genus $g$ and $r$ holes. The genus $g$ of each of these surfaces can be determined using the well-known formula for the Euler characteristic of such a surface, \begin{displaymath} \chi(S) = -2g -r+2, \end{displaymath} and solving the equation for the Thurston norm, \begin{displaymath} \| \phi \|_T = | \chi(S)| =d(2g+r-2), \end{displaymath} for $g$ using the value of $\| \phi \|_T$ from Equation \eqref{E:TN}. \end{rem} 
\section{The coincidence of Thurston and Alexander norms for graph links }
We now prove the main result of this article relating the Alexander and Thurston norms for graph links in homology 3-spheres. 
\begin{theorem}\label{T:TA} The Alexander and Thurston norms coincide for all irreducible graph links with two or more components.  \end{theorem}

\begin{rem}\label{R:P} By Corollary 8.3 of Eisenbud and Neumann \cite{Eis}, there is a unique minimal splice diagram for each graph link such that every edge weight is non-negative.  Hence in this proof without loss of generality we can assume that the edge weights of all the arrowhead and boundary vertices are non-negative; $\alpha_i \geq 0 , i=1, \dots ,r$ and $i=r+p+1, \dots , r+p+q$. Even more, since a boundary vertex with an edge weight equal to one represents a non-singular fiber, which is not an exceptional Seifert fiber, we can assume that $\alpha_i \geq 2$ for boundary vertices.\end{rem}

We prove the theorem by directly calculating the A-norm of the graph link $\textbf{L}$ using Equation \eqref{E:Al1} for the Alexander polynomial and comparing the result to Equation \eqref{E:TN} for the T-norm as given in \cite{Eis}.  We derive three new results which are used in the proof. Each of these results involves a sum over the $p$ splice components of the graph of $\textbf{L}$;\begin{enumerate}\item a \textit{T-norm decomposition formula}, 
\item  the Newton polyhedron of the Alexander polynomial of $\mathbf{L}$ as a Minkowski sum, and 
\item an \textit{A-norm decomposition formula}. \end{enumerate}

\subsection{T-norm decomposition formula} 
To show that the T-norm is a sum of the T-norms of each node of the graph of a graph link we need the following lemma.
\begin{lem}\label{L:BV} The linking numbers $l_{ji}$ of the arrowhead vertices into a node and into a boundary vertex attached to the same node of a graph link differ only by a factor of $\alpha_i$, where $\alpha_i$ is the weight of the edge in the graph connecting the boundary vertex indexed by $i$ to the node.\end{lem}

\begin{proof} The proof follows directly from the formula for the linking numbers given by Equation \eqref{E:ln}.  It says that to find the linking number between an arrowhead vertex and either a node or a boundary vertex we follow the path on the graph connecting the arrowhead vertex to the node or the boundary vertex. Along the way we multiply the product of all the signs of the nodes and also multiply all the edge weights of edges going into each node but not on the path.  Since a boundary vertex is attached to the node by hypothesis, it is clear that the paths between any arrowhead vertex and either the boundary vertex or the node it is attached to are the same except the path to the boundary vertex contains the edge connecting that vertex to the node.  Hence the weight of that edge, since it lies on the path, does not appear in the linking number between the arrowhead vertex and boundary vertex but it does in the linking number of that arrowhead vertex with the node.  Hence the two linking numbers differ exactly by the factor $\alpha_i$ which is the weight attached to the edge of boundary vertex $v_i$ indexed by $i$. \end{proof}
This lemma implies the following Corollary \ref{C:BV1} to Theorem \ref{T:TN} which gives a \textit{T-norm decomposition formula} for the irreducible graph link $\textbf{L}$:
\begin{cor}\label{C:BV1} The T-norm of graph link $\textbf{L}$, can be expressed as a sum of the T-norms of the $p$ splice components of the graph of $\textbf{L}$.
\begin{eqnarray}\label{E:TNBV} \| \phi \|_T & = &\sum_{i=r+1}^{r+p} (\tilde{\delta}_i - 2)\big{|} \sum_{j=1}^r \phi_j l_{ji}\big{|}\\ \nonumber &  =&\sum_{i=1}^p \| \phi \|_T^i, \end{eqnarray}  
where $\tilde{\delta}_i -2= \delta_i - 2- \sum_{k=1}^{q_i} \frac{1}{\alpha_k^i} > 0$ and
 $\| \phi \|_T^i$ is the contribution to the T-norm $\| \phi \|_T$ of the link from the $i$th splice component of the graph of the link. \end{cor}
\begin{rem}\label{R:1}We use the notation that $q_i$ denotes number of boundary vertices attached to the $i$th node, so that the boundary vertices $v_k^i$ attached to the $i$th node can be indexed by $k=1, \dots,q_i$.  In addition, we denote by  $\alpha_k^i$ the edge weight of the $k$th boundary vertex attached to the $i$th node.  In \cite{Eis}, the boundary vertices and the edge weights are ordered in a manner that does not specify which node the boundary vertex is attached to. Hence the notation $\alpha_i$ for an edge weight is used in \cite{Eis}, which we have refined to $\alpha_k^i$.\end{rem}

\begin{proof}  The expression for the T-norm of this corollary can be obtained from the equation of the T-norm given in \cite{Eis}, Equation \eqref{E:TN}, by a direct application of the Lemma \ref{L:BV} relating the linking numbers between arrowhead vertices to a node and the boundary vertices attached to the node.  In effect, the terms $\delta_i$ of the original expression are replaced by the terms $\tilde{\delta_i}$ and the sum over both nodes and boundary vertices is replaced by a sum over nodes only. To be precise we shall present this proof in detail.

Without loss of generality, we can prove the formula for the T-norm by showing it is true for the $(r+1)$th vertex which is also the first node that has by hypothesis $q_1$ boundary vertices attached to it. The equation given in \cite{Eis} for the T-norm of an irreducible graph link, Equation \eqref{E:TN}, can be written \begin{equation} \nonumber \| \phi \|_T = \sum_{i=r+1}^{r+p} (\delta_i - 2)\big{|} \sum_{j=1}^r \phi_j l_{ji} \big{|}  - \sum_{i=r+p+1}^{r+p+q} \big{|} \sum_{j=1}^r\phi_j l_{ji}\big{|} . \end{equation} We've split the sum into contributions from the nodes first and then the boundary vertices.  The contribution to the Thurston norm from the first node is \begin{equation} \nonumber \| \phi \|_T^1 = (\delta_{r+1} - 2)\big{|}\sum_{j=1}^r \phi_j l_{j(r+1)} \big{|}  - \sum_{i=r+p+1}^{r+p+q_1} \big{|}\sum_{j=1}^r \phi_j l_{ji} \big{|}. \end{equation} Applying the Lemma \ref{L:BV} that relates the linking numbers of the boundary vertices to the node we have that \begin{eqnarray} \nonumber  & &l_{ji} = \frac{1}{\alpha_i^{r+1}} l_{j(r+1)},\\ \nonumber &&\mbox{for }  j =1, \dots ,r  \mbox{ and } i=r+p+1, \dots ,r+p+q_1.  \end{eqnarray} We can combine the contribution from the node with the contributions from its boundary vertices to obtain \begin{eqnarray}\nonumber   \| \phi \|_T^1 &= &(\delta_{r+1} - 2- \sum_{k=r+p+1}^{r+p+q_1}\frac{1}{\alpha_k^{r+1}})\big{|} \sum_{j=1}^r \phi_j l_{j(r+1)} \big{|} \\ \nonumber
& = & (\tilde{\delta}_{r+1} - 2)\big{|}\sum_{j=1}^r \phi_j l_{j(r+1)} \big{|}.\end{eqnarray}  

By repeating this procedure indicated for the $(r+1)$th node on all $p$ of the nodes we obtain the Equation \eqref{E:TNBV} as claimed.

In order to show that $\tilde{\delta}_i - 2 > 0, \forall i$, we proceed by induction on the number of boundary vertices attached to the node. Without loss of generality we can show that this relation is true for the first node in order to prove it is true for all $p$ nodes. First, we use that $\delta_i \geq 3$ in the definition of a node given in \cite{Eis} and proceed by induction. Hence consider the first node indexed as $i=r+1$ and assume that it has only a single boundary vertex attached to it so that $q_1=1$. We have that $\delta_{r+1} \geq 3$ implies that \begin{equation}\nonumber \tilde{\delta}_{r+1} -2  = \delta_{r+1} - 2 - \frac{1}{\alpha_1^{r+1}} \geq 1 - \frac{1}{\alpha_1^{r+1}} > 0. \end{equation} Next let us assume that $\tilde{\delta}_{r+1}-2 > 0$ for $q_1 = n$ and we will show that this implies $\tilde{\delta}_{r+1}-2 > 0$ for $q_1 = n+1$. Adding the $(n+1)$th boundary vertex to the $(r+1)$th vertex, which is the first node, increases $\delta_{r+1}$ by one because of the additional edge into the node and also adds the term $(\alpha_{n+1}^{r+1})^{-1}$ to $\tilde{\delta}_i$. Adding this additional boundary vertex adds the strictly positive term $\frac{\alpha_{n+1}^{r+1} - 1}{\alpha_{n+1}^{r+1}} > 0 $ to $\tilde{\delta}_{r+1}$. By the inductive hypothesis $\tilde{\delta}_{r+1} -2 > 0$ without this additional boundary vertex. We find that we must have $\tilde{\delta}_{r+1} - 2 > 0$ after the addition of the strictly positive  contribution of the $(n+1)$th boundary vertex.  Hence by induction $\tilde{\delta}_{r+1} \geq 0$ for all values of $q_1$.\end{proof}

\subsection{Newton polyhedron of the Alexander polynomial of a graph link}
We also have to use another application of Lemma \ref{L:BV} in order to prove our fundamental theorem. It involves obtaining an expression for the Newton polyhedron $\Ne(\triangle^{\textbf{L}})$ of the graph link $\textbf{L}$. 

The \textit{Newton polyhedron} $\Ne(f)$ of a polynomial $f$ is the convex hull of the exponents of the Alexander polynomial $\triangle$.\footnote{The term \textit{Newton polytope} can also be used for the Newton polyhedron} Two Newton polyhedra can be added together using Minkowski addition: The \textit{Minkowski sum} of two polyhedra $K$ and $L$ is the set of all vector sums $x + y$ with $x\mbox{  } \in \mbox{  }K$ and $y \mbox{  } \in \mbox{  } L$.

The following proposition and corollary state two well-known properties of the Minkowski sum. The second corollary, although somewhat trivial, is new and will be essential in our proof of Theorem \ref{C:NP}. 
\begin{prop}(Gelfand, Kapranov and Zelevinsky \cite{Gel}, Prop.~6.1.2(b))\label{P:MS} The Newton polyhedron $\Ne(f\cdot f')$ of the product of two polynomials, $f$ and $f'$, is given by
\begin{equation}\nonumber \Ne(f\cdot f') = \Ne(f) + \Ne(f'). \end{equation}\end{prop}
Applying this proposition inductively to the product of $f$ with itself we obtain the following corollary.
\begin{cor}\label{CO:N} The Newton polyhedron $\Ne(f^n)$ of a polynomial $f$ to a power $n \in \mathbb{N}$, $f^n$, is given by
\begin{equation}\nonumber \Ne(f^n) = n \cdot \Ne(f).\end{equation}\end{cor}
\begin{cor}\label{C:FAC} Assume that $g$ divides $f$ for the rational polynomial $f/g$. Then the Newton polyhedron $\Ne(f/g)$ of this rational polynomial satisfies the relation \begin{equation}\label{E:Frac} \Ne(f) = \Ne(f/g) + \Ne(g). \end{equation}\end{cor}
\begin{proof} We use that $f$ can be written as a product of polynomials so that $f = \frac{f}{g} \cdot g$ and apply Proposition \ref{P:MS}.  \end{proof}

\begin{theorem}\label{C:NP} The Newton polyhedron $\Ne(\triangle^\textbf{L}(t_1, \dots ,t_r))$ of the graph link $\textbf{L}$, is given by \begin{eqnarray} \label{E:NP}  \Ne(\triangle^\textbf{L}(t_1, \dots ,t_r)) & = & \sum_{i=r+1}^{r+p} (\delta_i - 2 -  \sum_{j=1}^{q_i}\alpha_j^i) \Ne(t_1^{l_{1i}} \dots t_r^{l_{ri}}-1)\\ \nonumber
& =&\sum_{i=r+1}^{r+p} (\tilde{\delta}_i - 2 ) \Ne(t_1^{l_{1i}} \dots t_r^{l_{ri}}-1).\end{eqnarray}\end{theorem}
\begin{proof} By Corollary \ref{C:FAC}, given polynomials $f$ and $g$ such that $g$ divides $f$, their Newton polyhedra satisfy the equation $\Ne(f) = \Ne(f/g) + \Ne(g)$. If we set $f= \prod_{i=r+1}^{r+p} (t_1^{l_{1i}} \dots t_r^{l_{ri}}-1)^{\delta_i - 2}$ and $g = \prod_{i=r+p+1}^{r+p+q}(t_1^{l_{1i}} \dots t_r^{l_{ri}}-1)$, then $f/g$ is the Alexander polynomial $\triangle^{\textbf{L}}$ of $\textbf{L}$ as given in Equation \eqref{E:Al1}. We can prove this corollary by showing directly that $\Ne(f/g)= \sum_{i=r+1}^{r+p} (\delta_i - 2 -  \sum_{j=1}^{q_i}\alpha_j^i) \Ne(t_1^{l_{1i}} \dots t_r^{l_{ri}}-1)$ is a solution of the equation $\Ne(f) = \Ne(f/g) + \Ne(g)$. 

Substituting $\Ne(f)$ and $\Ne(g)$ into Equation \eqref{E:Frac} we obtain \begin{equation}\label{E:fg} \sum_{i=r+1}^{r+p} (\delta_i - 2) \Ne(t_1^{l_{1i}} \dots t_r^{l_{ri}}-1) = \Ne(f/g) + \sum_{i=r+p+1}^{r+p+q} \Ne(t_1^{l_{1i}} \dots t_r^{l_{ri}}-1).\end{equation} In this equation, we've also used that $\Ne((t_1^{l_{1i}} \dots t_r^{l_{ri}}-1)^{\delta_i-2}) = $\\$(\delta_i-2) \Ne(t_1^{l_{1i}} \dots t_r^{l_{ri}}-1)$  since $\Ne(f^n) = n \Ne(f), \forall n \in \mathbb{N}$, by Corollary \ref{CO:N}.  By Lemma \ref{L:BV}, the linking numbers of the boundary vertices are the same as the node they are attached to up to a factor of $(\alpha_j^i)^{-1}$. Using this result in $\Ne(g)$, we obtain \begin{equation} \nonumber \Ne(g) = \sum_{i=r+1}^{r+p} \sum_{j=1}^{q_i}\Ne \Bigl(t_1^{\frac{l_{1i}}{\alpha_j^i}} \dots t_r^{\frac{l_{ri}}{\alpha_j^i}}-1 \Bigr). \end{equation} The Newton polyhedra of this sum are all line segments with endpoints $(\frac{l_{1i}}{\alpha_j^i}, \dots , \frac{l_{ri}}{\alpha_j^i})$ and $(0, \dots, 0)$.  It is geometrically clear that $\alpha_j^i$ is a scaling factor of each of these segments which reduces the length of the segment but does not change its direction.  Hence we have the equality \begin{equation} \nonumber \Ne \Bigl(t_1^{\frac{l_{1i}}{\alpha_j^i}} \dots t_r^{\frac{l_{ri}}{\alpha_j^i}}-1 \Bigr) = \frac{1}{\alpha_j^i} \Ne(t_1^{l_{1i}} \dots t_r^{l_{ri}}-1). \end{equation} We now substitute this expression for $\Ne(g)$ along with our candidate solution for $\Ne(f/g)$ into Equation \eqref{E:fg} and obtain \begin{eqnarray}\nonumber  \sum_{i=r+1}^{r+p} (\delta_i - 2) \Ne(t_1^{l_{1i}} \dots t_r^{l_{ri}})& = & \sum_{i=r+1}^{r+p} (\delta_i - 2- \sum_{j=1}^{q_i} \frac{1}{\alpha_j^i}) \Ne(t_1^{l_{1i}} \dots t_r^{l_{ri}}-1)\\ \nonumber
&&\hspace{2mm} + \sum_{i=r+1}^{r+p}\sum_{j=1}^{q_i} \frac{1}{\alpha_j^i} \Ne(t_1^{l_{1i}} \dots t_r^{l_{ri}}-1)\\ \nonumber
&= &\sum_{i=r+1}^{r+p}\Bigl( \lambda_i \Ne(t_1^{l_{1i}} \dots t_r^{l_{ri}}-1)+ \lambda'_i \Ne(t_1^{l_{1i}} \dots t_r^{l_{ri}}-1)\Bigr) .\end{eqnarray} 
In this equation we've introduced the constants $\lambda_i = \delta_i - 2- \sum_{j=1}^{q_i} \frac{1}{\alpha_j^i}$ and $\lambda'_i = \sum_{j=1}^{q_i} \frac{1}{\alpha_j^i}$ each of which multiplies the same line segment $\Ne(t_1^{l_{1i}} \dots t_r^{l_{ri}}-1)$. The two terms on the right of this equation can be combined using Minkowski addition provided that $\lambda_i$ and $\lambda'_i$ are non-negative for all $i$. By Corollary \ref{C:BV1}, $\lambda_i  >0, \forall i$. As mentioned in Remark \ref{R:P}, we may assume, without loss of generality, that edge weights for boundary vertices satisfy the inequality $\alpha_j^i \geq 2, \forall i,j$. This implies that $\lambda'_i > 0, \forall i$. Hence since all the constants $\lambda_i$ and $\lambda'_i$ which multiply the same line segment are positive, we can use Minkowski addition to combine the two expressions on the right of this equation and the sums involving $\alpha_j^i$ cancel each other. 
\begin{eqnarray}\nonumber  \sum_{i=r+1}^{r+p} (\delta_i - 2) \Ne(t_1^{l_{1i}} \dots t_r^{l_{ri}}-1) &=&
 \sum_{i=r+1}^{r+p}( \lambda_i + \lambda'_i) \Ne(t_1^{l_{1i}} \dots t_r^{l_{ri}}-1)\\ \nonumber
&= &\sum_{i=r+11}^{r+p} \Bigl( \delta_i - 2- \sum_{j=1}^{q_i} \frac{1}{\alpha_j^i}+\sum_{j=1}^{q_i} \frac{1}{\alpha_j^i}\Bigr) \Ne(t_1^{l_{1i}} \dots t_r^{l_{ri}}-1)\\ \nonumber
& = &\sum_{i=r+1}^{r+p} (\delta_i - 2) \Ne(t_1^{l_{1i}} \dots t_r^{l_{ri}}-1).\end{eqnarray} Thus we see that the equality of Equation \eqref{E:fg} is satisfied.\end{proof}
A polyhedron which is the Minkowski sum of line segments is called a \textit{zonotope}.
\begin{cor}\label{C:z}The Newton polyhedron $\Ne(\triangle^\textbf{L}(t_1, \dots ,t_r))$ of the graph link $\textbf{L}$ is a zonotope consisting of $p$ line segments, one for each splice component ot the graph of $\textbf{L}$. \end{cor} 
\begin{proof} By Equation \eqref{E:NP}, $\Ne(\triangle^\textbf{L}(t_1, \dots ,t_r))$ is a Minkowski sum of the $p$ Newton polyhedra, $(\delta_i - 2) \Ne(t_1^{l_{1i}} \dots t_r^{l_{ri}}-1), i=1, \dots p$. The $i$th component of this sum, which is the Newton polyhedron of the $i$th splice component of the graph of $\textbf{L}$, has only two vertices, the endpoints $(\delta_i-2)(l_{1i}, \dots ,l_{ri})$ and $(0, \dots ,0)$. The convex hull of these two vertices is a line segment. \end{proof} 
\begin{rem} If the graph link $\textbf{L}$ has only one node, it is a Seifert link. By Corollary \ref{C:z} for $\Ne(\triangle^\textbf{L})$ with $p=1$, the Newton polyhedron of a Seifert link is a zonotope consisting of a single line segment. Dimca, Papadima and Suciu have found the same result in \cite{Dim} for the Newton polyhedron of a Seifert link by finding a coordinate system for which the Alexander polynomial is a function of only one variable, which they call its essential variable.
\end{rem}
\subsection{A-norm decomposition formula}  Since the Alexander norm is determined by the vertices of the Newton polyhedron $\Ne(\triangle)$, it can be viewed not only as a norm determined by an Alexander polynomial $\triangle$, but also as a norm determined by the Newton polyhedron $\Ne(\triangle)$ of the Alexander polynomial $\triangle$. This means that the notations \begin{eqnarray} \nonumber \| \cdot \|_A:&=& \| \cdot \|_{\triangle}\\ \nonumber &=&\| \cdot \|_{\Ne(\triangle)} \end{eqnarray} can also be useful.  In Long \cite{Lon}, we determine that the A-norm for the Newton polyhedron $\Ne(\triangle)$ is equal to the width function $w$ of $\Ne(\triangle)$ so that $\| \cdot \|_{\Ne(\triangle)} = w(\Ne(\triangle), \cdot)$. In the theory of polyhedra, the width function $w$ is a well-known function that can be proved to be \textit{Minkowski linear}. This means that for arbitrarily chosen polyhedra $P$ and $Q$ and all $\lambda \in \mathbb{R}^+$, it satisfies following:\begin{enumerate}
\item{Minkowski additivity : } $w(P+Q,\phi) =w(P,\phi) + w(Q, \phi)$.
\item{Minkowski scaling: } $w(\lambda P), \phi) = \lambda w(P, \phi) $. 
\end{enumerate} 
A proof of the Minkowski linearity of the width function can be found in Long \cite{Lont}.

The equivalence of the width function of $\Ne(\triangle)$ and the A-norm for $\Ne(\triangle)$ along with the Minkowski linearity of the width function imply that the A-norm has the following properties:
\begin{enumerate}
  \item  $\| \phi \|_{(\N(f) + \N(f'))} = \| \phi\|_{\N(f)} + \| \phi \|_{\N(f')}.$
  \item $\| \phi \|_{\lambda \N(f)} = \lambda \|\phi \|_{\N(f)}, \forall \lambda \in \mathbb{R}^+.$  
\end{enumerate}

By inductively applying Minkowski linearity, we obtain the following \textit{A-norm decomposition formula}.\begin{prop}\label{P:Andf} Assume that the Newton polyhedron $\Ne(\triangle)$ of the Alexander polynomial $\triangle$ can be written as a Minkowski sum of $k$ component Newton polyhedra, $\lambda_1 \Ne(f_1), \dots ,\lambda_k \Ne(f_k)$, with $\lambda_i \in \mathbb{R}^+, \forall i$, so that \begin{displaymath} \Ne(\triangle) = \sum_{i=1}^k \lambda_i \Ne(f_i). \end{displaymath} Then the A-norm for $\Ne(\triangle)$ is the sum of the A-norms of each component Newton polyhedron:\begin{equation} \label{E:ANDC}\| \phi \|_A =\| \phi \|_{\Ne(\triangle)}= \sum_{i=1}^k \lambda_i \| \phi \|_{\Ne(f_i)}. \end{equation} \end{prop}
\begin{rem} In Long \cite{Lon}, we derive a decomposition formula of the Alexander norm $\| \cdot \|_A=\| \cdot \|_{\triangle}$ for the polynomial $\triangle$ expressed as a sum involving the Alexander norms for each the irreducible factors of $\triangle$; $\| \phi\|_{\triangle} = \sum_{i=1}^k n_i \| \phi \|_{f_i}$ for $\triangle= f_1^{n_1} \dots f_k^{n_k}$ and $n_i   \in \mathbb{N},\forall i$.\end{rem}
By substituting the expression for the Newton polyhedron $\Ne(\triangle^{\textbf{L}})$ of $\textbf{L}$, Equation \eqref{E:NP}, into the A-norm decomposition formula, Equation \eqref{E:ANDC}, we obtain as a corollary the\textit{ A-norm decomposition formula for graph links}.
\begin{cor}\label{C:ANGL}The A-norm of graph link $\textbf{L}$, can be expressed as a sum of the A-norms of the $p$ splice components of the graph of $\textbf{L}$.
\begin{eqnarray}\label{E:ANGL}  \| \phi \|_A &= &\| \phi \|_{\Ne(\triangle^\textbf{L}(t_1, \dots ,t_r))} \\ \nonumber & = &\sum_{i=r+1}^{r+p}   (\tilde{\delta}_i - 2)\| \phi \|_{\Ne(t_1^{l_{1i}}  \dots t_r^{l_{ri}} -1)} \\ \nonumber & =&\sum_{i=1}^p \| \phi \|^i_A  \end{eqnarray} where $\| \phi \|_A^i$  denotes the contribution to the Alexander norm from the $i$th splice component of the graph of $\textbf{L}$.\end{cor}
\subsection{Proof of the main theorem}\begin{proof} 
We prove the theorem by showing that Equation \eqref{E:ANGL}, the A-norm decomposition formula for graph links, is equal to Equation \eqref{E:TNBV}, the T-norm decomposition formula. 
Each term in the sum of the A-norm decomposition formula for graph links involves the A-norm of a Newton polyhedron $\Ne(t_1^{l_{1i}}  \dots t_r^{l_{ri}} -1)$ which is a line segment with endpoints at the origin and at $(l_{1i}, \dots ,l_{ri})$. The A-norm for the line segment  $\Ne(t_1^{l_{1i}}  \dots t_r^{l_{ri}} -1)$  can  be obtained by direct application of the Definition \ref{D:fn} of the A-norm for the polynomial $t_1^{l_{1i}}  \dots t_r^{l_{ri}} -1$. This polynomial has two exponents $\alpha=(l_{1i}, \dots ,l_{ri})$ and $\beta=(0, \dots ,0)$.  The supremum which occurs in this definition must take its value on the difference of these two exponents since this polynomial has only two terms. Hence, we have that \begin{eqnarray}\nonumber \| \phi \|_A&=& \| \phi \|_{\Ne(t_1^{l_{1i}}  \dots t_r^{l_{ri}} -1)} = \| \phi \|_{(t_1^{l_{1i}}  \dots t_r^{l_{ri}} -1)}\\ \nonumber &= &|\phi(\alpha - \beta)| = |\phi(\alpha)| \\ \nonumber & =& \big{|}\sum_{j=1}^r \phi_j l_{ji} \big{|}.  \end{eqnarray}
 
After substituting this result that $\| \phi \|_{\Ne(t_1^{l_{1i}}  \dots t_r^{l_{ri}} -1)} = |\sum_{j=1}^r \phi_j l_{ji}|  $ into the Equation \eqref{E:ANGL} for the A-norm, we find that it is given by \begin{equation}\nonumber \|\phi \|_A =  \sum_{i=r+1}^{r+p} (\tilde{\delta}_i - 2) \big{|} \sum_{j=1}^r \phi_j  l_{ji}\big{|}  \end{equation} which is the same as the expression for the T-norm given in Equation \eqref{E:TNBV} by the T-norm decomposition formula.  \end{proof}

\subsection{The main theorem for graph knots}\begin{cor}\label{C:TAk} The Alexander and Thurston norms satisfy \begin{equation} \nonumber \label{E:knot} \| \phi \|_A - \| \phi \|_T = |\phi| \end{equation} for all irreducible graph knots. \end{cor}
\begin{proof} All the steps in the proof are the same as for an irreducible graph link with two or more components. However, since the knot has an extra factor of $(t-1)$ in its Alexander polynomial, which does not occur in a link with two or more components, this additional factor adds a summand of $\| \phi\|_{\Ne(t-1)}= | \phi |$ to the Alexander norm $\| \phi \|_A$ so that the two norms differ exactly by the amount $| \phi |$ as claimed; $\| \phi \|_A - \| \phi \|_T = | \phi|$.\end{proof}

\section{The equivalence of Alexander and Thurston norms for alternating links}
It has recently been proved by Ozsv\'ath and Szab\'{o} in \cite{Osz} that alternating links have identical Thurston and Alexander norms. We shall show in this section that the main result of this article, our Theorem \ref{T:TA}, that the two norms are equivalent for graph links complements their theorem on alternating links in its domain of applicability. 
\subsection{Theorem of Ozsv\'ath and Szab\'o}
The exact theorem of Ozsv\'ath and Szab\'o  in \cite{Osz} on alternating links is as follows.
\begin{theorem}[Ozsv\'ath and Szab\'o \cite{Osz}, Thm.\ 1.2, Cor.\ 1.3]\label{T:Ozs} Let $L \subset S^3$ be a link with $r$ components which admits a connected alternating projection.  The Newton polytope $\Ne(\triangle^L)$ of the Alexander polynomial $\triangle^L$ of $L$ is equivalent up to a scale factor of two to the dual Thurston unit norm polyhedron $\mathcal{B}_T^*$ of the complement of $L$. \end{theorem}

This theorem is proved using Floer homology. From it we can directly deduce the following as a corollary.

\begin{cor} The Thurston and Alexander norms coincide for links $L$ as in the above theorem. \end{cor}
\begin{proof}   In Long \cite{Lon}, it is shown that in the coordinate system with the Newton polytope $\Ne(\triangle^L)$ centered at the origin, the dual Newton polytope $\Ne(\triangle^L)^*$ is equivalent to the A-norm unit ball $\mathcal{B}_A$ up to a scale factor of two; $\Ne(\triangle^L)^* = 2 \mathcal{B}_A$. By this Theorem \ref{T:Ozs}, the  Newton polytope is  equal to the dual T-norm unit ball up to a scale factor of two; $\Ne(\triangle^L)  = 2 \mathcal{B}_T^*$. This equivalence is preserved under the duality transformation so that $\Ne(\triangle^L)^*  = 2 \mathcal{B}_T$ also. Hence these two results combined imply that the A-norm unit ball is identical to the T-norm unit ball; $\mathcal{B}_A = \mathcal{B}_T$. This equivalence implies the two norms must coincide for every cohomology class $\phi$ so that $\| \phi \|_A = \| \phi \|_T$ as claimed.
\end{proof}

\subsection{Menasco-Thurston theorem}
An alternating link has a projection where the overcrossings and undercrossings alternate. Using results of Thurston \cite{Thu2}, Menasco showed that alternating links are either of torus type or else the link complement has a hyperbolic geometry in \cite{Men}. Hence these results can be summarized in the following \textit{Menasco-Thurston theorem on alternating links}.

\begin{theorem}[Menasco, \cite{Men}, Cor. 2] \label{T:Men}If L is a non-split prime alternating link which is not a torus link, then $S^3 \setminus L$ has a complete hyperbolic structure (of finite volume).\end{theorem}
\begin{rem} It is well-known that not only are the $(2,n)$-torus links the only alternating torus links, but also that their link complements do not have a hyperbolic geometry. \end{rem}

This theorem has important applications because of the theory of the JSJ decomposition of 3-manifolds discussed in the next section.  
 \section{Coincidence of Thurston and Alexander norms and the JSJ decomposition}
We conclude by the Menasco-Thurston Theorem \ref{T:Men} that our Theorem \ref{T:TA} is original in that it applies to a significant set of links for which the Ozsv\'ath-Szab\'{o} Theorem \ref{T:Ozs} does not, all the graph links which are built up by splicing Seifert-fibered components along tori except the $(2,n)$-torus links. We also conclude that our theorem on the equivalence of Alexander and Thurston norm for graph links complements the Ozsv\'ath-Szab\'{o} Theorem \ref{T:Ozs} on the equivalence of these norms for alternating links with respect to the JSJ decomposition.  According to the JSJ decomposition for links, the link complement can be decomposed into two types of pieces; Seifert and hyperbolic links. A hyperbolic link is defined as follows:
\begin{defn} A link $L$ is hyperbolic if its complement $X$ admits a complete hyperbolic structure of finite volume.\footnote{In \cite{Eis} on page 24, the terminology \textit{simple} is used in place of hyperbolic: ``A simple link is defined to be an irreducible link $\textbf{L}= (\Sigma,L)$ with the property that any incompressible torus in $\Sigma\setminus \interior N(L)$ ($\interior N(L)$ denotes the interior of a closed tubular neighborhood of the link $L$ in the oriented homology 3-sphere $\Sigma$) is boundary parallel. By Thurston's hyperbolization theorem \cite{Thu2}, this is equivalent to saying that that the link complement $\Sigma\setminus L$ admits a complete hyperbolic structure of finite volume, except possibly if $L=\phi$ and $\Sigma$ is not sufficiently large.'' } \end{defn}

Having defined the two types of pieces, Seifert-fibered and hyperbolic, we can now define the JSJ decomposition for the complement of a link in $S^3$ as follows.
\begin{defn}[Eisenbud and Neumann \cite{Eis}, pg. 25]Let $L$ be a link in $S^3$. Then the link complement $X = S^3 \setminus L$ admits an essentially unique splice decomposition along tori into irreducible pieces which are either Seifert-fibered or hyperbolic called the \textit{JSJ decomposition}.\end{defn}
A piece of the decomposition is said to be irreducible if and only if it does not have a splice decomposition into component pieces which are Seifert-fibered or hyperbolic.

Our goal is to combine our main Theorem \ref{T:TA} together with the Ozsv\'ath-Szab\'{o} Theorem \ref{T:Ozs} and the Menasco-Thurston Theorem \ref{T:Men} to deduce a more general theorem in the context of the JSJ decomposition of the link exterior. In order to do this we define the splice notation of the JSJ decomposition and use it to define the notion of a JSJ component. 
\begin{defn}\label{D:spl1}The notation \textit{$\mathcal{N}_i \cup_T \mathcal{N}_j$} denotes the \textit{splice} of two components $\mathcal{N}_i$ and $\mathcal{N}_j$ indexed by $i$ and $j$ of the JSJ decomposition of a link across a torus $T$.\end{defn}

Our generalized theorem states that the Thurston and Alexander norms are the same for links in $S^3$ which are either graph links or alternating links.
\begin{theorem} Let $L$ be a link in $S^3$.  Let the JSJ decomposition of the link complement $S^3 \setminus L$ consist of $k$ irreducible pieces $\mathcal{N}_i$ so that \begin{equation} \nonumber S^3 \setminus L = \mathcal{N}_1 \cup_T \mathcal{N}_2 \cup_T \dots \cup_T \mathcal{N}_k. \end{equation}
Then the Thurston and Alexander norms of $L$ coincide if either of the following two conditions are met:
\begin{enumerate}
\item $L$ is a graph link so that all the $\mathcal{N}_i$ are Seifert-fibered.
\item $L$ is an alternating link which is not a $(2,n)$-torus link, so that all the $\mathcal{N}_i$ are hyperbolic. \end{enumerate}\end{theorem}
\begin{proof} The proof of the result implied by the first condition can be found in the proof of our main Theorem \ref{T:TA} that the Thurston and Alexander norms are the same for graph links in a homology sphere. The 3-sphere $S^3$ is a homology 3-sphere and graph links are precisely those links for which all the irreducible pieces of the JSJ decomposition of its complement in $S^3$ are Seifert-fibered. We can conclude that Theorem \ref{T:TA} implies that for all the links in $S^3$ whose JSJ decomposition consists solely of Seifert-fibered pieces, the Alexander and Thurston norms coincide as claimed.

The proof of the result implied by the second condition follows directly from the Ozsv\'ath-Szab\'{o} Theorem \ref{T:Ozs} which states that the Thurston and Alexander norms coincide for alternating links. Further by the Menasco-Thurston Theorem \ref{T:Men} these links, except for the $(2,n)$-torus links, can have only hyperbolic pieces in their JSJ decomposition since the link complement must have a hyperbolic geometry. \end{proof}

\section{Characteristic hyperplanes of the Thurston norm unit ball}
In this section, we introduce the reduced T-norm unit ball which can be decomposed into two sets; the fibered and non-fibered cohomology classes.  Next, we show that for graph links the set of non-fibered cohomology classes form reduced characteristic hyperplanes through the origin which pass through all the faces of the reduced T-norm unit ball except the top-dimensional faces. 
\subsection{Thurston and Alexander norm unit balls}\label{S:Ub}  We summarize here some of the results from Long \cite{Lon} on the A-norm unit ball. In general, for a 3-manifold $M$ with Betti number $b_1$,  the Alexander norm unit ball $\mathcal{B}_A$  is not a bounded set.  The dimension $b_e$ of the Newton polyhedron $\Ne(\triangle)$, with $b_e \leq b_1$, determines the number of coordinates we call \textit{essential coordinates}. These coordinates determine a decomposition of $\mathbb{R}^{b_1}$ into the product of two vector subspaces, the vector space spanned by the essential coordinates $V_e = \mathbb{R}^{b_e}$ and the space spanned by the non-essential coordinates $V_{ne} = \mathbb{R}^{b_1-b_e}$; $\mathbb{R}^{b_1} = V_e \times V_{ne}$. We call $V_e$ the \textit{essential vector space} and $V_{ne}$ the \textit{non-essential vector space}.  For the irreducible graph link $\textbf{L}$ with $r$ link components, the first Betti number is equal to the number of link components; $b_1=r$. 

In a coordinate system of $\mathbb{R}^{b_1}$ with such a vector space decomposition, the A-norm is a function of only the essential coordinates.  When restricted to essential coordinates, the A-norm unit ball $\mathcal{B}_A$ becomes the reduced A-norm unit ball $\mathcal{\tilde{B}}_A$ which is a convex polyhedron in $V_e$. The A-norm unit ball $\mathcal{B}_A$ is the product set of   $\mathcal{\tilde{B}}_A$ and $V_{ne}$;  \begin{equation} \label{E:RN} \mathcal{B}_A = \mathcal{\tilde{B}}_A \times V_{ne} =\mathcal{\tilde{B}}_A \times \mathbb{R}^{b_1-b_e}. \end{equation} This relation implies that  $\tilde{\mathcal{B}}_A$ is the projection of $\mathcal{B}_A$ onto the vector space $V_e$; $\tilde{\mathcal{B}}_A = \pi_{V_e} \mathcal{B}_A$ where $\pi_{V_e}$ denotes the projection mapping of $\mathbb{R}^r$ onto $V_e$. Since the A-norm and T-norm coincide for graph links, we introduce the\textit{ reduced T-norm unit ball $\tilde{\mathcal{B}}_T$} which is equal to the reduced A-norm unit ball for these links; $\tilde{\mathcal{B}}_T=\tilde{\mathcal{B}}_A$ . 

\subsection{3-manifold fibrations}
The set of cohomology classes of the T-norm unit ball $\mathcal{B}_T$ can be broken up into two subsets; fibered and non-fibered.  An integer-valued cohomology class is said to be fibered or non-fibered depending on whether or not the surface $S$ it represents is the fiber $F$ of a fibration of the 3-manifold $M$ over the circle $S^1$. We define precisely what it means to say a cohomology class $\phi \in H^1(M;\mathbb{R})$ represents a fibration in the following definition.
\begin{defn}We call a real-valued cohomology class $\phi \in H^1(M;\mathbb{R})$ a \textit{fibered cohomology class} if every integer-valued class on the ray from the origin through the class $\phi$ represents a surface $S$ through Poincar\'e duality which is the fiber $F$ of a fibration of $M$ over the circle $S^1$; $F \to M \to S^1$.  Otherwise, $\phi$ is said to be a \textit{non-fibered cohomology class}.\end{defn}
\begin{rem} The fiber which the integer-valued class $\phi$ represents need not be connected. In fact, if the greatest common divisor of the components of a class $\phi$ is equal to $d$, then $\phi$ represents a fiber which is the disjoint union of $d$ identical connected components.\end{rem}The set of fibered classes consists of a set of open cones called \textit{Thurston cones}. These cones have their cone point at the origin and pass through the interior of a subset of the top-dimensional faces of the T-norm unit ball. A top-dimensional face of a polyhedron is called a facet. 
\begin{defn} A facet $\mathcal{F}$ of the reduced T-norm unit ball $\mathcal{\tilde{B}}_T$ is a \textit{fibered facet} if every non-zero cohomology class $\phi$ in the cone $\mathcal{C}_{\mathcal{F}}$ pointed at the origin and passing through the interior of $\mathcal{F}$ is a fibered cohomology class. \end{defn}
A fibered facet determines a Thurston cone as follows.
\begin{defn} The cone $\mathcal{C}_{\mathcal{F}}$ pointed at the origin of the reduced T-norm unit ball $\mathcal{\tilde{B}}_T$ through the facet $\mathcal{F}$ is a \textit{Thurston cone} if the facet $\mathcal{F}$ is a fibered facet. \end{defn}The following \textit{Thurston cone theorem} divides the facets of the reduced T-norm unit ball into two subsets; fibered and non-fibered facets. Each fibered facet determines a Thurston cone so that the total set of fibered cohomology classes $\textbf{F}$ is the union of the Thurston cones, minus the origin.
\begin{theorem}[Thurston \cite{Thu}, Thm.\ 5]\label{T:Thu}  Let $M$ be a compact oriented $3$-manifold.  The set $\textbf{F}$ of cohomology classes representable by non-singular (non-vanishing) closed one-forms is some union of the cones on open (top-dimensional) faces of $\mathcal{B}_T$, minus the origin.  The set of elements in $H^1(M;\mathbb{Z})$ whose Lefschetz (Poincar\'e) dual is represented by a fiber of a fibration consists of all lattice points in $\textbf{F}$.\end{theorem}
\begin{rem} The above statement of the Thurston cone theorem from \cite{Thu} assumes that the T-norm is a norm rather than a semi-norm so that the unit ball $\mathcal{B}_T$ is a convex polyhedron.  In the case that the T-norm is not a norm, so that there are non-zero classes   with zero T-norm, we use the reduced T-norm unit ball $\mathcal{\tilde{B}}_T$ in this theorem instead since this unit ball is a convex polyhedron.\end{rem}
The following \textit{Tischler fibration theorem} on the fibrations of $n$-dimensional manifolds over the circle explains the use of closed, non-vanishing one-forms in the Thurston cone theorem in determining the set of fibered classes.
\begin{theorem}[Tischler \cite{Tis}, Theorem 1]\label{T:Tis} Let $M^n$ be a closed $n$-dimensional manifold.  Suppose $M^n$ admits a closed, non-vanishing one-form.  Then $M^n$ is a fiber bundle over $S^1$.\end{theorem}
\begin{rem} In the proof of this theorem, a one-form $\omega$ is taken to be an element of the first deRham cohomology class of $M^n$; $\omega \in H_{DR}^1(M;\mathbb{R})$. By assuming that $\omega$ is a closed, non-vanishing one-form with rational coefficients, a fibration of $M^n$ over $S^1$ is constructed. Since the rational-valued cohomology classes are dense in $H^1(M;\mathbb{R})$, there is always a rational-valued class $\epsilon$-close to any real-valued cohomology class for any $\epsilon > 0$. \end{rem}
\begin{cor}  If $\omega \in H^1(M^n;\mathbb{R})$ is a closed, non-vanishing one-form with rational coefficients, then $\omega$ induces a fibration map, $\Lambda :M^n \to S^1$, of the fiber bundle of the above theorem over the circle.\end{cor} 

\subsection{The fibration criterion for graph links}The following \textit{Eisenbud-Neumann fibration criterion} gives necessary and sufficient conditions for a cohomology class of a graph link to be fibered. 
\begin{theorem}[Eisenbud and Neumann \cite{Eis}, Theorem 11.2] A necessary and sufficient condition for the cohomology class $\phi$ of the irreducible graph link $\textbf{L}$ to be fibered is that \begin{equation}\nonumber  \sum_{j=1}^{r} \phi_j l_{ji}\neq 0 \mbox{ for } i=r+1, \dots ,r+p+q\end{equation}\end{theorem}

By applying our Lemma \ref{L:BV} that the linking number $l_{ji}$ between the $j$th arrowhead vertex and the $i$th vertex, which is either a node or a boundary vertex, differ by only a constant when the boundary vertex is attached to the node, we can reduce the set of conditions of the above theorem for the graph link $\textbf{L}$ to be fibered from $p+q$ (number of nodes and boundary vertices) to $p$ (number of nodes only).
\begin{cor} A necessary and sufficient condition for the cohomology class $\phi$ of the irreducible graph link $\textbf{L}$ to be fibered is that \begin{equation} \nonumber \sum_{j=1}^{r} \phi_j l_{ji}\neq 0 \mbox{ for } i=r+1, \dots ,r+p \end{equation}\end{cor}
\begin{proof} By Lemma \ref{L:BV}, $q$ of the $p+q$ conditions given in the theorem are redundant. Hence the number of conditions can be reduced from $p+q$ to $p$.\end{proof}  
As a second corollary, we can deduce that the set of non-fibered cohomology classes is given by $p$ hyperplanes which we call the \textit{characteristic hyperplanes}.
\begin{cor} \label{T:SC}A cohomology class $\phi$ of the irreducible graph link $\textbf{L}$ is non-fibered if and only if it is contained in one of $p$ hyperplanes $\mathcal{H}_i, i=1, \dots ,p,$ through the origin in $\mathbb{R}^r$.  There is one hyperplane for each node of the graph and it is given explicitly by \begin{displaymath}\mathcal{H}_i = \big{\lbrace} \phi \in H^1(X;\mathbb{R}) \mid \sum_{j=1}^r \phi_j l_{ji} = 0 \big{\rbrace}.\end{displaymath}\end{cor}

\begin{proof}The statement used in this corollary for a cohomology class to be non-fibered, $\sum_{j=1}^r \phi_j l_{ji} = 0$ for at least one of the nodes, can be obtained as the negation of the statement of the previous corollary for a cohomology class to be fibered.  The T-norm for the $r$ component graph link $\textbf{L}$ is defined in a space of dimension $r$, so the condition $\sum_{j=1}^r \phi_j l_{ji} = 0$ determines an $r-1$ dimensional hyperplane through the origin in this space as claimed. Since there are $p$ such conditions there must be $p$ hyperplanes. \end{proof}
We can characterize the characteristic hyperplanes $\mathcal{H}_i$ in a slightly different manner using our T-norm decomposition formula given by Equation \eqref{E:TNBV}.
\begin{cor} The characteristic hyperplanes $\mathcal{H}_i, i=1, \dots, p,$ of the above corollary are the sets \begin{equation} \label{E:Ch} \mathcal{H}_i = \lbrace \phi \in H^1(X;\mathbb{R}) \mid \| \phi\|_T^i =0\rbrace. \end{equation}\end{cor}
\begin{proof} By Equation \eqref{E:TNBV}, the T-norm $\| \phi\|_T$ of the cohomology class $\phi$ can be written as a sum of the T-norms of each of the splice components of the graph of the link $\textbf{L}$, $\| \phi \|_T^i, i=1, \dots ,p$, where $\| \phi\|_T^i= (\tilde{\delta}-2)\sum_{j=1}^r \phi_j l_{ji}$. Since $(\tilde{\delta}-2)$ is always positive, it follows that $\| \phi\|_T^i=0$  if and only if $\sum_{j=1}^r \phi_j l_{ji}=0$. \end{proof}
\subsection{The reduced characteristic hyperplanes}   In this section, we use a coordinate system for $\mathbb{R}^r$ which decomposes as a product space so that $\mathbb{R}^r = V_e \times V_{ne}$ as discussed in Section \ref{S:Ub}. We then determine what happens to the characteristic hyperplanes when they are projected onto the vector space $V_e$ spanned by the essential coordinates. 
\begin{prop}\label{P:Tne}  The T-norm $\| \phi \|_T$ of the class $\phi$ for the irreducible graph link $\textbf{L}$ is a function of only the $b_e$ essential components of $\phi$. \end{prop}
\begin{proof} As mentioned in Section \ref{S:Ub}, it is shown in Long \cite{Lon} that the A-norm is a function of only the $b_e$ essential coordinates which span the essential vector space $V_e$. By our main theorem, Theorem \ref{T:TA}, the A-norm coincides with the T-norm for the irreducible graph link $\textbf{L}$ so the T-norm must also be a function of only the $b_e$ essential coordinates for this link.  \end{proof}
\begin{defn}The projection of the $i$th characteristic hyperplane $\mathcal{H}_i$  onto the essential vector space $V_e$ is the $i$th \textit{reduced characteristic hyperplane $\tilde{\mathcal{H}}_i$}:
\begin{equation}\label{E:rch}  \tilde{\mathcal{H}}_i  = \pi_{V_e} \mathcal{H}_i  =\lbrace \tilde{\phi} \in V_e \mid \| \phi \|_T^i =0\rbrace. \end{equation} \end{defn}
\begin{prop}\label{P:Eh} The $i$th reduced characteristic hyperplane is isomorphic as a vector space to the quotient space of the $i$th characteristic hyperplane by the non-essential vector space.\begin{displaymath}  \tilde{\mathcal{H}}_i \cong  \mathcal{H}_i/ V_{ne}. \end{displaymath}\end{prop}
\begin{proof} The defining condition $\| \phi \|_T^i=0$ of $\mathcal{H}_i$ does not depend on the non-essential components of $\phi$ by Proposition \ref{P:Tne}. \end{proof}
\begin{rem} We can use the type of same argument to show the reduced T-norm unit ball is also isomorphic to a quotient space. The defining condition of $\mathcal{\tilde{B}}_T$ that $\| \phi \|_T \leq 1$ does not depend on the non-essential coordinates by Proposition \ref{P:Tne}. This implies that $\tilde{\mathcal{B}}_T \cong \mathcal{B}_T/V_{ne}$.\end{rem} 
\begin{lem}\label{L:Nf} The cohomology class $\phi \in \mathbb{R}^r$ of the graph link $\textbf{L}$ is non-fibered if and only if its projection, $\tilde{\phi} = \pi_{V_e} \phi$, onto $V_e$ lies in one of the $p$ reduced characteristic hyperplanes $\tilde{\mathcal{H}}_i, i=1, \dots ,p$. \end{lem}
\begin{proof}  This is a direct consequence of the equivalence relation of Proposition \ref{P:Eh} and that each non-fibered class $\phi \in \mathbb{R}^r$ lies in one of the characteristic hyperplanes $\mathcal{H}_i,i=1, \dots ,p$, by Corollary \ref{T:SC}.\end{proof}    

\subsection{The fibered facets of the reduced Thurston norm unit ball} We use the reduced characteristic hyperplanes and the Thurston cone theorem to prove the following new theorem on graph link fibrations.
\begin{theorem} \label{T:Ff}  Every facet $\mathcal{F}$ of the reduced T-norm unit ball $\mathcal{\tilde{B}}_T$ for the irreducible graph link $\textbf{L}$ is a fibered facet. \end{theorem}
\begin{proof}We can assume, as discussed in Section \ref{S:Ub}, that the reduced T-norm unit ball $\mathcal{\tilde{B}}_T$ is a convex polyhedron of dimension $b_e\leq r$ in the essential vector space $V_e = \mathbb{R}^{b_e}$. By Lemma \ref{L:Nf}, the set of non-fibered cohomology classes is made up of the union of the $p$ reduced characteristic hyperplanes $\tilde{\mathcal{H}}_i$ in $V_e$.  By their definition in Equation \eqref{E:rch}, each of the these $b_e-1$ dimensional hyperplanes pass through the origin. On the other hand, none of the facets of  $\mathcal{\tilde{B}}_T$ intersect the origin since it must contain the origin and the facets lie on its boundary. Hence, if an $b_e-1$ dimensional hyperplane $\tilde{\mathcal{H}}_i$ intersects an $b_e-1$ dimensional facet $\mathcal{F}$, the intersection, which is the set $\tilde{\mathcal{H}}_i \cap \mathcal{F}$, must be a set of at most $b_e-2$ dimensions. Thus it is not possible for a facet of $\mathcal{\tilde{B}}_T$ to have its interior consist solely of non-fibered cohomology classes. Since each facet must consist solely of either fibered or non-fibered classes by the Thurston cone theorem, we can conclude every facet is a fibered facet as claimed.\end{proof}

\section{The sample graph link $\textbf{L}_{EN}$}
\begin{figure}[h]
\includegraphics[width=9cm,height = 6cm]{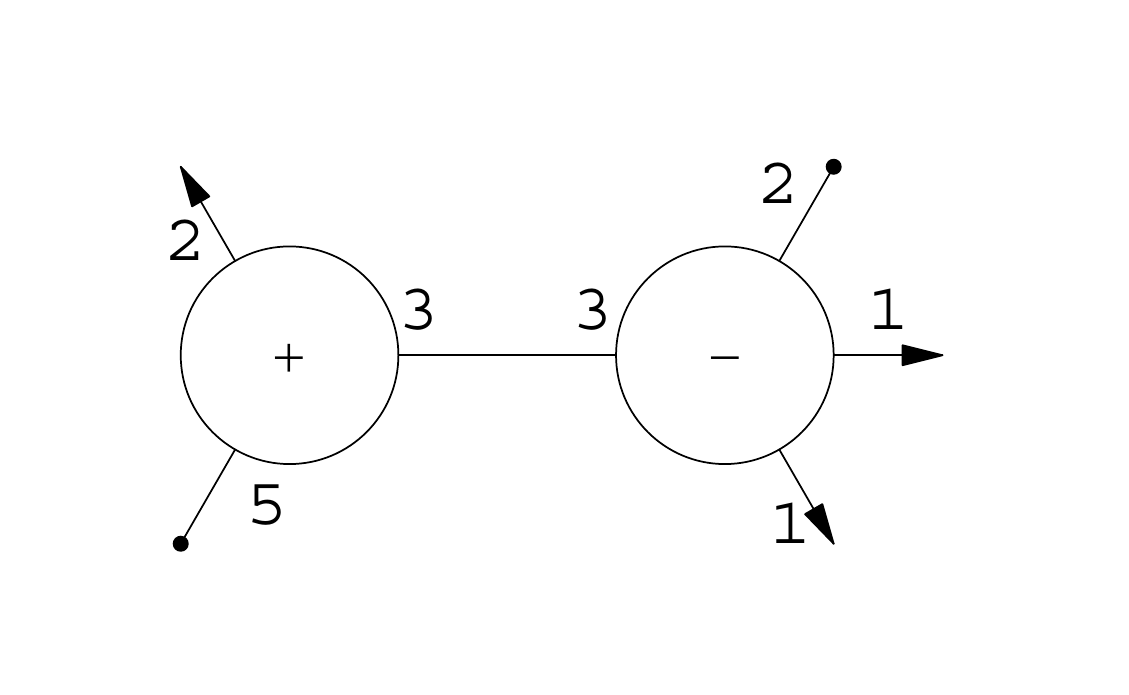}
\caption{The splice diagram of the graph link $\textbf{L}_{EN}$ from page 68 of Eisenbud and Neumann \cite{Eis}. This figure shows the splice of link $\textbf{L}=(\Sigma(+;2,3,5), S_1 \cup S_2)$ on the left with the link $\textbf{L}'=(\Sigma(-;1,1,3,2),S_1' \cup S_2' \cup S_3')$ on the right. The splice is from $S_2$ of $\textbf{L}$ to $S_3'$ of $\textbf{L}'$ and results in  the composite link $\textbf{L}_{EN}$. }\label{F:S2}
\end{figure}

We refer to the graph link whose splice diagram is shown in Figure \ref{F:S2} as $\textbf{L}_{EN}$ and use it as a standard example for calculations.   
\subsection {The Alexander polynomial for graph link $\textbf{L}_{EN}$} \label{Ex:LEN}
By using the splice diagram for $\textbf{L}_{EN}$ shown in Figure \ref{F:S2}, the algorithm of Equation \eqref{E:ln} to find the linking numbers $l_{ji}$ and Equation \eqref{E:Al1} for the Alexander polynomial of a graph link, we obtain the following expression for the Alexander polynomial of $ \textbf{L}_{EN}$.
\begin{equation} \label{E:GL1} \triangle^{\textbf{L}_{EN}}(t_1,t_2,t_3) = \frac{(t_1^{15} t_2^{-20}t_3^{-20} - 1)(t_1^{-10} t_2^{-6} t_3^{-6} - 1)^2} {( t_1^3 t_2^{-4} t_3^{-4} - 1)(t_1^{-5} t_2^{-3} t_3^{-3}-1)}. \end{equation} 
The factored form of this polynomial is
\begin{eqnarray}\nonumber \label{E:GL2}  \triangle^{\textbf{L}_{EN}}(t_1,t_2,t_3)& = & -t_1^{-15} t_2^{-25} t_3^{-25} (t_1^5 t_2^3 t_3^3 - 1)(t_1^5 t_2^3 t_3^3 + 1)^2\\ \nonumber
 &&(t_1^{12} + t_1^9 t_2^4 t_3^4 + t_1^6 t_2^8 t_3^8 + t_1^3t_2^{12}t_3^{12} + t_2^{16} t_3^{16}).  \end{eqnarray} 
Since all Alexander polynomials are defined only up to multiplication by monomial factors, the monomial factor with the negative exponents can be dropped from $\triangle^{\textbf{L}_{EN}}$. The Newton polyhedron $\Ne(\triangle^{\textbf{L}_{EN}})$ for this link is the convex hull of the exponents. Although the exponents lie in a three-dimensional space, the convex hull of the exponent set is only two-dimensional so that $b_1=3$ and $b_e=2$. The variables $t_2$ and $t_3$ appear in each monomial term of $\triangle^{\textbf{L}_{EN}}$ with the same exponent.  Hence a monomial isomorphism given by $\tilde{t}_1 = t_1, \tilde{t}_2= t_2 t_3 $ and $\tilde{t}_3 = t_3$ eliminates the variable $\tilde{t}_3$ from the set of exponents. In the new coordinate system, the variables $\tilde{t}_1$ and $\tilde{t}_2$ are essential variables and $\tilde{t}_3$ is the only non-essential variable. After doing this, we find that $\Ne(\triangle^{\textbf{L}_{EN}})$ in essential variables is the convex hull of the set of vertices $\lbrace(0,16),(12,0),(27,9),(15,25) \rbrace $. 
\subsection{The T-norm and its unit ball for graph link $\textbf{L}_{EN}$}\label{EX:GL}
We use the T-norm decomposition formula, Equation \eqref{E:TNBV}, to determine the T-norm of $\textbf{L}_{EN}$.  From the Alexander polynomial $\triangle^{\textbf{L}_{EN}}$ given in Equation \eqref{E:GL1}, we find that the linking numbers $l_{ji}$ of the three arrowhead vertices, indexed by $j=1,2,3$, into the first node, indexed by $i=4$, on the left in Figure \ref{F:S2} are $l_{14}=15$ and $l_{24}=l_{34}=-20$ and into the second node, indexed by $i=5$,  are $l_{15}= -10$ and $l_{25}=l_{35}= -6$. By counting the number of edges into the two nodes, we find that $\delta_4= 3$ and $\delta_5=4$. The edge weights of the single boundary vertex into each node are $\alpha_1^4= 5$ and $\alpha_1^5=2$. Thus $\tilde{\delta}_4 = \delta _4 - \alpha_1^4 = 3-1/5=14/5$ and $\tilde{\delta_5} = \delta_5 - \alpha_1^5= 4 - 1/2=7/2$. Substituting these quantities into Equation \eqref{E:TNBV}, we find that \begin{eqnarray}\nonumber \label{E:TNGL}   \|\phi\|_T  & = &\frac{4}{5}|15\phi_1 -20 \phi_2 -20 \phi_3| + \frac{3}{2}|-10 \phi_1 -6 \phi_2-6 \phi_3| \\ \nonumber & =&4|3\phi_1-4(\phi_2+\phi_3)|+3 |5\phi_1+3(\phi_2 + \phi_3)|\\ \nonumber &= &\| \phi \|_T^1 + \| \phi \|_T^2.  \end{eqnarray}

The unit ball $\mathcal{B}_T$  is the following set of points: \begin{equation} \label{E:Rtn1} \nonumber \mathcal{B}_T = \lbrace \phi \in \mathbb{R}^3  \mid \| \phi \|_T = 4|3\phi_1-4(\phi_2+\phi_3)|+3 |5\phi_1+3(\phi_2 + \phi_3)| \leq 1 \rbrace \end{equation}  The change of coordinates from $(t_1, t_2,t_3)$ to $(\tilde{t}_1, \tilde{t}_2, \tilde{t}_3)$ given above induces the change of coordinates in cohomology from $(\phi_1, \phi_2, \phi_3)$ to $(\tilde{\phi}_1, \tilde{\phi}_2,\tilde{\phi}_3)$ given by $\tilde{\phi}_1 = \phi_1$,  $\tilde{\phi}_2 = \phi_2 + \phi_3$ and $\tilde{\phi}_3 = \phi_3$. 
In the new coordinate system, the Thurston norm is 
\begin{equation} \nonumber \| \phi \|_T= 4| 3 \tilde{\phi}_1 - 4 \tilde{\phi}_2| + 3 | 5 \tilde{\phi}_1 + 3 \tilde{\phi}_2|.\end{equation}
The coordinates of the reduced Thurston norm unit ball $\tilde{\mathcal{B}}_T$ are the essential coordinates $\tilde{\phi}_1$ and $\tilde{\phi}_2$ and they span the essential vector space $V_e = \mathbb{R}^2$;    \begin{equation} \label{E:Rtn2}  \tilde{\mathcal{B}}_T = \lbrace  \tilde{\phi} \in \mathbb{R}^2 \mid \| \phi \|_T= 4| 3 \tilde{\phi}_1 - 4 \tilde{\phi}_2| + 3 | 5 \tilde{\phi}_1 + 3 \tilde{\phi}_2|  \leq 1 \rbrace \end{equation} By substituting $b_e=2$ and $b_1=3$ into Equation \eqref{E:RN} relating the A-norm unit ball to the reduced A-norm unit ball, and using that the A-norm is equal to the T-norm, we find that \begin{displaymath} \mathcal{B}_T = \tilde{\mathcal{B}}_T \times V_{ne} =\tilde{\mathcal{B}}_T \times \mathbb{R}. \end{displaymath}
 
The vertices of the reduced T-norm unit ball determined by this inequality can be found by setting each absolute value to zero in Equation \eqref{E:Rtn2} and solving the equality $\| \phi \|_T=1$. Since the T-norm unit ball is symmetric with respect to the origin, we also know that if $v$ is a vertex then $-v$ is also a vertex. \begin{enumerate} \item $4| 3 \tilde{\phi}_1 - 4 \tilde{\phi}_2| = 0$ and $3 |5 \tilde{\phi}_1 + 3 \tilde{\phi}_2| = 1 \Rightarrow (4/77,3/77)$ and $(-4/77,-3/77)$ are vertices.
\item $3|5 \tilde{\phi}_1 + 3 \tilde{\phi}_2| = 0$ and  $4| 3 \tilde{\phi}_1 - 4 \tilde{\phi}_2| = 1 \Rightarrow (-3/116,5/116)$ and  $(3/116,-5/116)$ are vertices.
\end{enumerate} The full set of vertices is \begin{equation}\nonumber  \ver(\tilde{\mathcal{B}}_T) 
=\lbrace (4/77,3/77),(-4/77,-3/77),(-3/116,5/116),(3/116,-5/116)\rbrace.  \end{equation} The convex hull of this set of vertices is shown in Figure \ref{F:GL1}. 

\subsection{The characteristic hyperplanes for graph link $\textbf{L}_{EN}$}
\begin{figure}[h]
\includegraphics[width=6cm,height = 6cm]{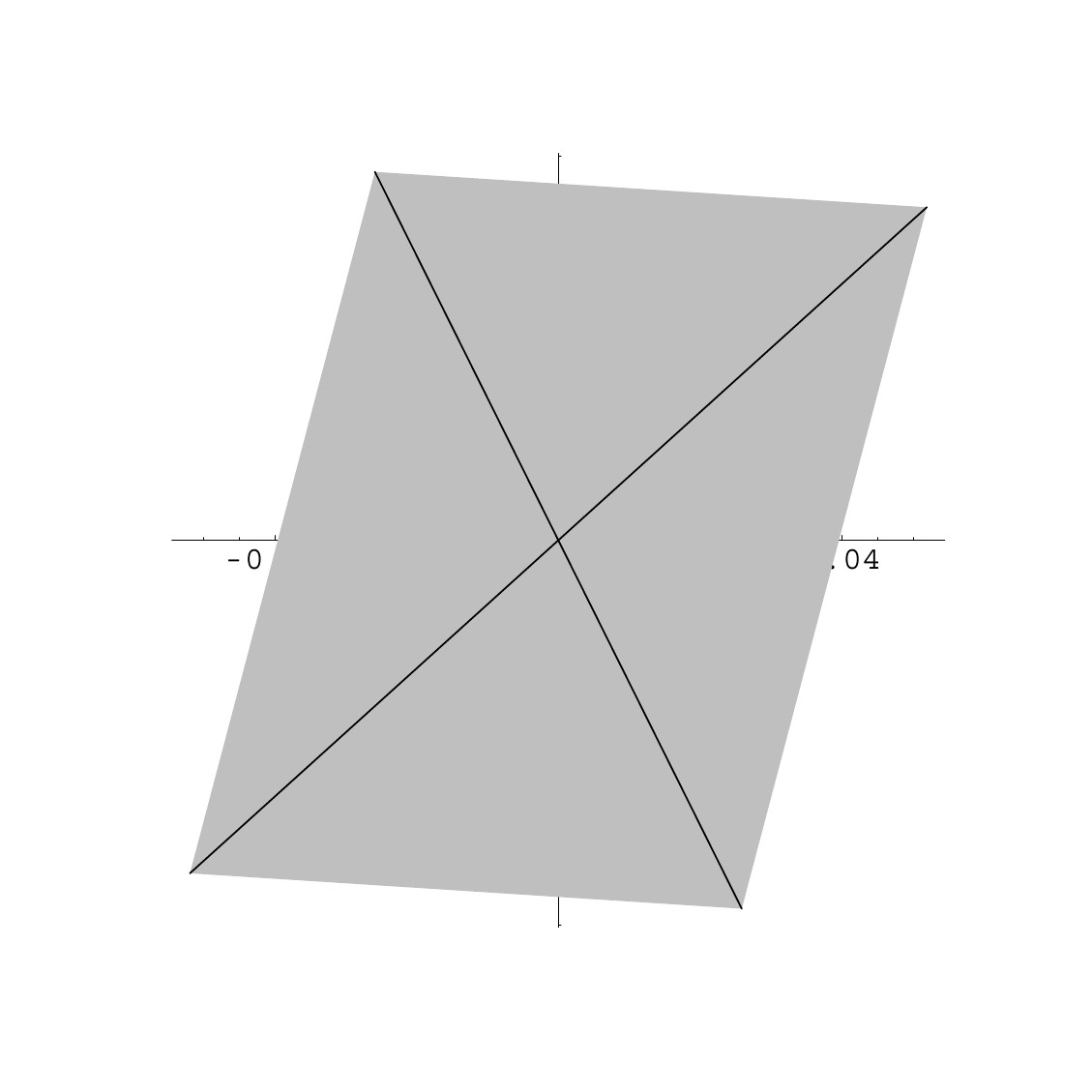}
\caption{The intersection of the reduced characteristic hyperplanes $\tilde{\mathcal{H}}_1$ and $\tilde{\mathcal{H}}_2$ and the reduced T-norm unit ball $\tilde{\mathcal{B}}_T$ of the link $\textbf{L}_{EN}$.}\label{F:GL1}
\end{figure}  
As an example of determining the non-fibered set of cohomology classes of the T-norm unit ball, we calculate this set for the sample graph link $\textbf{L}_{EN}$. The characteristic hyperplanes are found by setting the two Thurston norms for the splice components, $\| \phi \|_T^1$ and $\| \phi \|_T^2$, in the Equation \eqref{E:TNGL} for the Thurston norm equal to zero.  They are \begin{eqnarray} \nonumber  && \mathcal{H}_1 = \lbrace \phi \in \mathbb{R}^3 \mid 3\phi_1-4(\phi_2 +\phi_3) = 0 \rbrace \\ \nonumber
&& \mathcal{H}_2 = \lbrace \phi \in \mathbb{R}^3 \mid  5\phi_1+3(\phi_2 +\phi_3) = 0 \rbrace.  \end{eqnarray} There are two characteristic hyperplanes, $\mathcal{H}_1$ and $\mathcal{H}_2$, because there are two nodes in the graph shown in Figure \ref{F:S2}. Each of these is a two-dimensional plane through the origin. By projecting the characteristic hyperplanes $\mathcal{H}_1$ and $\mathcal{H}_2$ onto the essential vector space $V_e$, we obtain the reduced characteristic hyperplanes $\tilde{\mathcal{H}}_1$ and $\tilde{\mathcal{H}}_2$. They are lines determined by the following two equations: \begin{eqnarray} \nonumber  && \tilde{\mathcal{H}}_1 =\lbrace \tilde{\phi} \in \mathbb{R}^2 \mid  3 \tilde{\phi}_1 - 4 \tilde{\phi}_2 =0 \rbrace\\ \nonumber
&& \tilde{\mathcal{H}}_2 =\lbrace \tilde{\phi} \in \mathbb{R}^2 \mid 5 \tilde{\phi}_1 +3  \tilde{\phi}_2=0 \rbrace.  \end{eqnarray} The intersection of these two lines with the reduced T-norm unit ball $\tilde{\mathcal{B}}_T$  is the set of non-fibered cohomology classes of $\tilde{\mathcal{B}}_T$. As predicted by our Theorem \ref{T:Ff}, every facet of $\tilde{\mathcal{B}}_T$ is a fibered facet so there are four Thurston cones. Thus the Figure \ref{F:GL1} shows the two lines pass through all the vertices of $\tilde{\mathcal{B}}_T$. In accordance with the Thurston cone Theorem \ref{T:Thu}, since each of the four cones of $\tilde{\mathcal{B}}_T$ is a Thurston cone, these four vertices make up the total set of non-fibered faces; the faces with codimension greater than one of $\tilde{\mathcal{B}}_T$.

\textbf{Acknowledgements }\newline
This article is based on Ph.D. thesis research conducted at Northeastern University under the direction of Alexander~I.~Suciu. Further details can be found in Long \cite{Lont}, \cite{Lon}.

\end{document}